\documentclass[12pt]{amsart}
\usepackage[utf8]{inputenc}
\usepackage{amsmath,amssymb,amsthm, graphics}
\usepackage{amsthm}
\usepackage{diagmac2}
\usepackage{textcomp}
\usepackage[alphabetic]{amsrefs}
\usepackage[all,cmtip]{xy}
\usepackage{yfonts}
\usepackage{fullpage}
\usepackage{stmaryrd}
\usepackage{mathrsfs}
\usepackage{caption}

\newcommand{\bP}{\mathbb{P}}

\newcommand{\bF}{\mathbb{F}}

\newtheorem{theorem}{Theorem}[section]
\newtheorem{Def}[theorem]{Definition}
\newtheorem{Teo}[theorem]{Theorem}
\newtheorem*{Teo1}{Theorem 1.2}
\newtheorem*{Teo2}{Theorem 1.4}

\newtheorem{Cor}[theorem]{Corollary}

\newtheorem{Question}[theorem]{Question}
\newtheorem*{Question*}{Question 1}
\newtheorem*{Question2*}{Question 2}

\begin{document}

\title{Sharp Bertini Theorem for Plane Curves over Finite Fields}
\author{Shamil Asgarli}

\begin{abstract}
We prove that if $C$ is a reflexive smooth plane curve of degree $d$ defined over a finite field $\bF_q$ with $d\leq q+1$, then there is an $\bF_q$-line $L$ that intersects $C$ transversely. We also prove the same result for non-reflexive curves of degree $p+1$ and $2p+1$ where $q=p^{r}$. 
\end{abstract}

\maketitle

\section{Introduction}

A classical theorem of Bertini states that if $X$ is a smooth quasi-projective variety in $\bP^n$ defined over an \textit{infinite} field $k$, then a general hyperplane section of $X$ is smooth. Specializing to the case when $C\subseteq\bP^2$ is a smooth plane curve, it follows that there exists a line $L$ (defined over $k$) such that $L$ intersects $C$ transversely, meaning that $C\cap L$ consists of $d$ distinct geometric points where $d=\deg(C)$. But when $k=\bF_q$ is a finite field, it is possible to have a smooth plane curve $C\subseteq\bP^2$ such that every line $L$ defined over $\bF_q$ is tangent to the curve $C$ (see Example 2.A below). Moreover, Poonen's Bertini Theorem \cite{Poonen}*{Theorem 1.2} guarantees that such smooth curves, where all the $\bF_q$-lines are tangent, do exist in every sufficiently large degree (see Example 2.B below). With a view toward an effective version of Poonen's theorem, one can ask the following:

\begin{Question}\label{GeneralQuestion} Suppose $C \subseteq\bP^2$ is a smooth plane curve defined over $\bF_q$. Let $d=\deg(C)$. What conditions on $q$ and $d$ will ensure that there is a line $L\subseteq\bP^2$ defined over $\bF_q$ such that $L$ meets $C$ transversely?
\end{Question}

Let us call $L$ a \textbf{good line} if $L$ meets $C$ transversely. We expect that if $q$ is large with respect to $d$, then good lines will exist. Indeed, if $q\geq d(d-1)$, then the dual curve $C^{*}$ cannot be space-filling, i.e. $C^{*}(\bF_q) \neq (\bP^2)^{*}(\bF_q)$. This is because $\deg(C^{*})\leq d(d-1)\leq q$ and a curve of degree of at most $q$ cannot go through all the points of $(\bP^2)^{*}(\bF_q)$. Any point in $(\bP^2)^{*}(\bF_q)\setminus C^{\ast}(\bF_q)$ represents a good line $L\subseteq \bP^2$ defined over $\bF_q$. A generalization of this observation to higher dimensions is proved by Ballico \cite{Bal}*{Theorem 1}. 

In this paper, we improve the quadratic bound $q\geq d(d-1)$ to the linear bound $q\geq d-1$.

\begin{Teo}\label{Teo:reflexive} If $C$ is a smooth reflexive plane curve defined over $\bF_q$ with $\deg(C)\leq q+1$, then there is an $\bF_q$-line $L$ such that $L$ intersects $C$ transversely. \end{Teo}
The theorem is sharp in a sense that the statement cannot be improved to $q\geq d-2$. There is a counter-example when $q=d-2$ (see Example 2.A). The ``reflexive'' assumption on $C$ is same as saying that $C$ has finitely many flex points (see Section 2). As a natural follow-up, we may ask:

\begin{Question}\label{question:nonreflexive}
Does Theorem \ref{Teo:reflexive} hold when C is non-reflexive?
\end{Question}

We prove a partial result in this direction:

\begin{Teo}\label{Teo:nonreflexive} Let $C$ be a smooth non-reflexive plane curve of degree $p+1$ or $2p+1$ defined over $\bF_q$ where $q=p^{r}$ with $r\geq 2$. Then there is an $\bF_q$-line $L$ such that $L$ intersects $C$ transversely. 
\end{Teo}

Finally, in the last section of the paper (Section 4), we focus exclusively on Frobenius non-classical curves, which are non-reflexive curves of special kind. As we will see, Question \ref{question:nonreflexive} in this case is equivalent to a statement about collinear $\bF_q$-points on the curve. 

\textbf{Conventions.} In order to avoid various pathologies, we will assume throughout the paper that the characteristic of the field is $p>2$. 

\textbf{Acknowledgements.} I would like to thank my advisor Brendan Hassett for the unwavering support and constant encouragement. I am especially grateful to Felipe Voloch for suggesting to investigate the case of degree $7$ non-reflexive curves over $\bF_9$, which led to Theorem \ref{Teo:nonreflexive}. I also thank Dan Abramovich, Dori Bejleri, Herivelto Borges, Pol van Hoften and Giovanni Inchiostro for insightful discussions and useful comments on the manuscript.

\section{Reflexive Curves} 
In this section we review the theory of reflexive plane curves, and prove Theorem \ref{Teo:reflexive}. \\ 
If $C$ is a plane curve defined over a field $k$, we can consider the Gauss map $\varphi: C\to (\bP^2)^{\ast}$ that associates to each smooth point $p$ of $C$ its tangent line. The dual curve $C^{\ast}$ is defined to be the closure of $\varphi(C)$ inside $(\bP^2)^{\ast}$. By looking at the Gauss map for the dual curve, we get $\varphi': C^{\ast}\to C^{\ast\ast}$. In what follows, we will identify $\bP^{2}$ and $(\bP^{2})^{\ast\ast}$.

\begin{Def} The curve $C$ is called \textbf{reflexive} if $C=C^{\ast\ast}$ and $\varphi'\circ\varphi: C\to C^{\ast\ast}$ is the identity map.
\end{Def}

A theorem of Wallace \cite{Wallace} asserts that $C$ is reflexive if and only if $\varphi$ is separable. As a result, all smooth plane curves in characteristic zero are reflexive. Recall that a point $P$ of $C$ is called a \textbf{flex point} if the tangent line at $P$ meets the curve $C$ at $P$ with multiplicity at least $3$. When $\operatorname{char}(k)=p>2$, we have the following characterization: $C$ is reflexive if and only if $C$ has finitely many flex points \cite{Par}*{Proposition 1.5}.

Before we prove Theorem \ref{Teo:reflexive}, here are some counter-examples of smooth curves $C$ where all the lines defined over $\bF_q$ are tangent to $C$ (so that no good line exists). 

\textbf{Example 2.A.} Let $C$ be a smooth plane curve with $\deg(C)=q+2$ such that $\# C(\bF_q) = \#\bP^2(\bF_q)$. Such curves exist, and have been extensively studied by Homma and Kim \cite{HomKim}. For such a curve $C$, every $\bF_q$-line $L$ intersects $C$ at $q+2$ points (counted with multiplicity). But $q+1$ of these points are already accounted by the points of $L(\bF_q)=\bP^1(\bF_q)$. Thus, the residual intersection multiplicity results from $L$ being tangent to $C$ at one of the $\bF_q$-points. 

\textbf{Example 2.B.} Fix a finite field $\bF_q$. Let $\{L_1, ..., L_{q^2+q+1}\}$ be all the $\bF_q$-lines in the plane. Pick distinct (geometric) points $P_i\in L_i$ for each $i$. The condition that $C$ is tangent to $L_i$ at $P_i$ is a statement about vanishing of the first few coefficients in the Taylor expansion at these finitely many points. By applying Poonen's Bertini theorem with Taylor conditions \cite{Poonen}*{Theorem 1.2}, there exists some $d_0$ such that for \textit{every} $d\geq d_0$, there exists a smooth plane curve $C\subseteq \bP^2$ of degree $d$ such that $L_i$ is tangent to $C$ at $P_i$. In particular, all $\bF_q$-lines $L\subseteq\bP^2$ are tangent to $C$. A closer inspection of the proof reveals that the integer $d_0$ is in the order of $q^2$ (essentially because we imposed $q^2+q+1$ local conditions).

We will now prove the main theorem of the present paper.

\begin{Teo1} If $C$ is a smooth reflexive plane curve defined over $\bF_q$ with $\deg(C)\leq q+1$, then there is an $\bF_q$-line $L$ such that $L$ intersects $C$ transversely. \end{Teo1}

\begin{proof} Let $\Phi$ be the Frobenius map defined on points by $\Phi([X:Y:Z])=[X^q:Y^q:Z^q]$. We will write $T_P(C)$ for the tangent line to $C$ at a (geometric) point $P$. Set $$N = \# \{P\in C(\overline{\bF_q}): \Phi(P)\in T_P(C)\}$$ which is finite because $C$ is reflexive \cite{HefVol}. The following inequality is proved in \cite{HKT}*{Theorem 8.41}:
\[
2\cdot \#C(\bF_q) + N \leq d(q+d-1) \tag{$\ast$}
\]
under the assumption that $C$ has finitely many flex points and that characteristic of the field is $p>2$. This is the step where we use the hypothesis that $C$ is reflexive.

Assume, to the contrary, that every $\bF_q$-line is tangent to the curve $C$ at some (geometric) point. Let us divide these lines into two groups: if $L$ is tangent to $C$ at an $\bF_q$-rational point, we will call $L$ a \textbf{rational tangent}. Otherwise, we will call $L$ a \textbf{special tangent}. Since every $\bF_q$-line is tangent to $C$, and there are $q^2+q+1$ lines defined over $\bF_q$, we get 
$$
\# \{\text{rational tangents}\} + \# \{\text{special tangents}\}  = q^2 + q + 1 
$$
and
$$
\# \{\text{rational tangents}\} \leq \# C(\bF_q)
$$
Now, if $L$ is a special tangent, it is tangent to the curve $C$ at a non-$\bF_q$-point $P$. Then $L$ is also tangent to $C$ at $P, \Phi(P), \Phi^2(P), ..., \Phi^{e-1}(P)$ where $e=[k(P):\bF_q]$ is the degree of the point $P$. Since $e\geq 2$, the line $L$ contributes at least $2$ elements to $N$. As a result, 
$$
2\cdot \# \{\text{special tangents}\} \leq N
$$
Combining all the inequalities above, we obtain that
\begin{align*}
q^2 + q + 1 &= \# \{\text{rational tangents}\} + \# \{\text{special tangents}\} \\
&\leq \#C(\bF_q) + \frac{N}{2} \leq \frac{1}{2} d(q+d-1) \tag{using $(\ast)$} \\ 
&\leq \frac{1}{2} (q+1)(q+(q+1)-1) = \frac{1}{2}(q+1)(2q) = q^2+q
\end{align*}
which is a contradiction.
\end{proof}

When $q=p$ is a prime, every smooth curve of degree at most $p$ is reflexive. Moreover, Pardini \cite{Par}*{Proposition 3.7} has shown that every smooth non-reflexive curve of degree $p+1$ (over any field of characteristic $p$) is projectively equivalent to the curve given by the equation $xy^p + yz^p + zx^p=0$. For this curve, many good lines exist. For instance, take two $\bF_p$-points on the curve, and join them with a line $L$. Then $L$ will intersect $C$ transversely.

Consequently, we deduce the result for all smooth plane curves over $\bF_p$ where $p$ is prime.

\begin{Cor}\label{Fp} If $C$ is a smooth plane curve defined over $\bF_p$ with $\deg(C)\leq p+1$ where $p$ is a prime, then there is an $\bF_p$-line $L$ such that $L$ intersects $C$ transversely. 
\end{Cor}

\section{Non-reflexive curves}

In this section, we will restrict attention to non-reflexive curves and prove Theorem \ref{Teo:nonreflexive}. \\
Let $C\subseteq\bP^2 $ be a smooth non-reflexive curve defined over $\bF_q$ with $q=p^{r}$ where $r\geq 2$. Pardini \cite{Par}*{Corollary 2.4} has shown that $C$ is defined by an equation of the form:
$$
a^p x + b^p y + c^p z = 0
$$
where $a, b, c\in\bF_{q}[x,y,z]$ are homogeneous polynomials of degree $t\geq 1$. In particular, $\deg(C)=tp+1$. 

We establish a Bertini-type theorem for the case $t=1$ and $t=2$.

\begin{Teo2} Let $C$ is a smooth non-reflexive plane curve of degree $p+1$ or $2p+1$ defined over $\bF_q$ where $q=p^{r}$ with $r\geq 2$. Then there is an $\bF_q$-line $L$ such that $L$ intersects $C$ transversely. 
\end{Teo2}

\begin{proof} When $\deg(C)=p+1$, then $C$ is projectively equivalent to the curve given by the equation $xy^p + yz^p + zx^p=0$, for which many good lines $L$ exist (see the discussion before Corollary \ref{Fp}). For the rest of the proof, we will assume that $\deg(C)=2p+1$. Since $C$ is non-reflexive, by \cite{Par}*{Corollary 4.3} the degree of the dual curve is 
$$
\deg(C^{\ast}) = \frac{d(d-1)}{p} = \frac{(2p+1)(2p)}{p} = 4p+2
$$
For $p\geq 5$, we observe that $\deg(C^{\ast})=4p+2\leq p^2\leq q$, so $C^{\ast}$ cannot contain all of $(\bP^2)^{\ast}(\bF_q)$, and hence any point $L \in  (\bP^2)^{\ast}(\bF_q) \setminus C^{\ast}(\bF_q)$ will be a desired line that intersects $C$ transversely. 

When $p=3$, the inequality $\deg(C^{\ast})=4p+2=14\leq p^r=q$ still holds for $r\geq 3$. The only case that requires a separate analysis is $(p, r)=(3, 2)$, which corresponds to degree $2\cdot 3 + 1 = 7$ curve defined over $\bF_{3^2}=\bF_{9}$. The rest of the proof is devoted to studying this remaining case.

Let $C$ be a smooth non-reflexive curve of degree $7$ defined over $\bF_9$. Assume, to the contrary, that all the lines defined over $\bF_9$ are tangent to $C$. Following the same terminology used in the proof of Theorem \ref{Teo:reflexive}, we call $L$ a \textbf{rational tangent} if $L$ is tangent to $C$ at some $\bF_9$-point. Otherwise, $L$ is called a \textbf{special tangent}. Since $C$ is non-reflexive, each tangent line $L$ must intersect the curve at the tangency point with multiplicity $\geq 3$ (Proposition 1.5 in \cite{Par}). It follows that: 
\begin{enumerate} 
\item If $L$ is a rational tangent, then $L\cap C$ contains at most five $\bF_9$-points. 
\item If $L$ is a special tangent, then $L\cap C$ contains a conjugate pair of $\bF_{81}$-points and a single $\bF_9$-point. In symbols, $L\cap C=\{Q, Q^{\sigma}, P\}$ where $Q\in\bP^2(\bF_{81})\setminus\bP^2(\bF_9)$ and $P\in\bP^2(\bF_9)$. 
\end{enumerate}

Consider the following incidence correspondence of points and lines,
$$
\mathcal{I}= \{(P,L): L\in(\bP^2)^{\ast}(\bF_9) \text{ and } P\in (C\cap L)(\bF_9)  \}
$$
Each $P\in C(\bF_9)$ is contained in $q+1 = 10$ different $\bF_9$-lines. Therefore, $\#\mathcal{I} = \# C(\bF_9)\cdot 10$. On the other hand, using (1) and (2) above, each special tangent $L$ contributes $1$ point, while each rational tangent $L$ contributes at most $5$ points to $\#\mathcal{I}$. Thus, $\#\mathcal{I}\leq S + 5R$ where $S$ and $R$ are the number of special and rational tangents, respectively. We deduce that
$$
\# C(\bF_9)\cdot 10 \leq S +  5R 
$$
Since $\# C(\bF_9)\geq R$, we get $10 R \leq S + 5R$, which implies $5R\leq S$. Since $S+R =9^2+9+1 = 91$, we have $5(91-S)\leq S$, so that $S\geq \frac{5\cdot 91}{6} = 75.8333...$. Thus, $S\geq 76$.  

Next, take any rational tangent $L_0$. Every special tangent line intersects $L_0$ in one of its ten $\bF_9$-points. Since $\frac{S}{10}\geq \frac{76}{10} > 7$, there exists $P_0\in L_0(\bF_q)$ such that there are at least $8$ special tangent lines that pass through $P_0$. By looking at the ten $\bF_9$-lines passing through $P_0$, we can estimate $\#C(\bF_9)$ as follows. Each of the $8$ special tangents will contribute at most $1$ rational point, while the remaining (at most $2$) rational tangents will contribute at most $5$ rational points. Thus, one gets $\# C(\bF_9) \leq 8+2\cdot 5 = 18$. Consider the incidence correspondence:
$$
\mathcal{J} = \{(P,L): L \text{ is a special tangent and } P\in (C\cap L)(\bF_9)  \}
$$
By (1) above, every special tangent contains exactly one $\bF_9$-point of $C$, so that $\# \mathcal{J} = S$. As a result,
$$
S = \#\mathcal{J} = \sum_{P\in C(\bF_9)} \#\{\text{special tangents passing through } P\}
$$
Since
$$
\frac{S}{\#C(\bF_9)} \geq \frac{76}{18}>4
$$
there exists a point $P\in C(\bF_9)$ such that at least $5$ special tangents pass through $P$. Consider the corresponding line $P^{\ast}$ in the dual space $(\bP^2)^{\ast}$, which consists of all lines passing through $P$. Let us look at the intersection of the line $P^{\ast}$ and the dual curve $C^{\ast}$ inside $(\bP^2)^{\ast}$. The intersection has all the ten $\bF_9$-points of $P^{\ast}$ since all the $\bF_9$-lines are tangent to $C$. However, each of the special tangents is bitangent to $C$, so it is a node in $C^{\ast}$, and hence will contribute $2$ to the intersection. It follows that $P^{\ast}\cap C^{\ast}$ has at least $5\cdot 2 + 5 = 15$ intersections, contradicting the fact that $\deg(C^{\ast})=14$. 
\end{proof}

\textbf{Remark.} As we saw above, the hardest part of the proof is the case $p=3$. This answers a question of Felipe Voloch, who asked in a private communication, whether or not there exists a transverse line for a degree $7$ smooth non-reflexive curve defined over $\bF_{9}$. The small primes still persist when we try to extend Theorem \ref{question:nonreflexive} to non-reflexive curves of degree $3p+1$. Indeed, if $C$ is a smooth non-reflexive curve of degree $3p+1$, then  $\deg(C^{\ast}) = \frac{(3p+1)(3p)}{p} = 9p+3 \leq p^2\leq q$ for $p\geq 11$; the usual argument shows that $(C^{\ast})(\bF_q)\neq (\bP^2)^{\ast}(\bF_q)$, implying that good lines exist for $p\geq 11$. However, the main difficulty lies with the primes $p=3, 5, 7$.   

\section{Connection to Frobenius non-classical curves}

In this section, we observe the implications of a Bertini-type theorem for a special class of non-reflexive curves, known as Frobenius non-classical curves. 

\begin{Def}
Let $C\subseteq\bP^2 $ be a smooth plane curve defined over $\bF_q$. Then $C$ is called \textbf{Frobenius non-classical} if $\Phi(P)\in T_{P}(C)$ for every $P$, where $T_{P}(C)$ is the tangent line to $C$ at the point $P$, and $\Phi: \bP^2 \to \bP^2$ is the $q$-th power Frobenius map.  
\end{Def}

We should remark that the usual definition of Frobenius non-classical is stated differently (by looking at the order sequence of $C$), but the definition given above is equivalent in the case of smooth plane curves \cite{HefVol}*{Proposition 1}. 

\textbf{Example.}
Let $C$ be the curve defined over $\bF_{q^2}$ by the equation
$$
x^{q+1} + y^{q+1} + z^{q+1} = 0
$$
It can be checked that $C$ is a smooth Frobenius non-classical curve for $\bF_{q^2}$. 

If $C$ is a smooth Frobenius non-classical plane curve of degree $d$ defined over $\bF_{q}$ where $q=p^{r}$, then it is known that $C$ is non-reflexive \cite{HefVol}*{Proposition 1} and $\sqrt{q}+1 \leq d \leq \frac{q-1}{q'-1}$ where $q'$ is the generic order of contact of the curve with a tangent line \cite{HefVol}*{Propositions 5 and 6}. In particular, $\deg(C)\leq q-1$ always holds. So Question \ref{question:nonreflexive} is equivalent to:

\begin{Question}\label{Frob}
If $C$ is a smooth Frobenius non-classical plane curve defined over $\bF_q$, does there exist an $\bF_q$-line $L$ such that $L$ intersects $C$ transversely?
\end{Question}

The existence of such a line $L$ can be verified for the curve $x^{q+1} + y^{q+1} + z^{q+1} = 0$, and more generally, for the curve given by the equation
$$
x^{q^{n-1}+\cdots + q +1} + y^{q^{n-1}+\cdots + q +1} + z^{q^{n-1}+\cdots + q+1} = 0
$$
where $n\geq 2$. These curves are indeed smooth and Frobenius non-classical with respect to the field $\bF_{q^{n}}$ \cite{HefVol}*{Theorem 2}. 

If the Question \ref{Frob} has an affirmative answer, then it implies that there is a line $L$ defined over $\bF_q$ such that $L\cap C$ consists of $d=\deg(C)$ distinct $\bF_q$-rational points. Indeed, if $L$ contains a non-$\bF_q$-point $Q$, then we observe that $Q, \Phi(Q)\in T_{Q}(C)$ (since $C$ is Frobenius non-classical) and $Q, \Phi(Q)\in L$ (as $L$ is defined over $\bF_q$), implying that $L=T_{Q}(C)$ is a tangent line. Thus, any good (transverse) line $L$ intersects $C$ at $\deg(C)$ distinct $\bF_q$-points. This allows us to reformulate Question \ref{Frob} as follows:

\begin{Question}\label{Frob2}
If $C$ is a smooth Frobenius non-classical plane curve defined over $\bF_q$, then does $C$ have $d=\deg(C)$ many $\bF_q$-rational points on a line?
\end{Question}

The Question \ref{Frob2} is motivated by the fact that Frobenius non-classical curves have many $\bF_q$-points. In fact, the $\bF_q$-points on these curves have been used in \cite{GPTU} and \cite{Borges} to construct certain complete arcs in the plane. Moreover, the following theorem due to Hefez and Voloch \cite{HefVol}*{Theorem 1} gives the exact the number of $\bF_q$-points on \textit{any} smooth Frobenius non-classical plane curve:

\begin{Teo}\label{HF}(Hefez-Voloch)
If $C\subseteq\bP^2$ is a smooth Frobenius non-classical curve of degree $d$ defined over $\bF_q$, then
$$
\# C(\bF_q) = d(q-d+2)
$$
\end{Teo}

We can apply Theorem \ref{HF} directly to get an estimate on the number of collinear points of $C$. Consider the incidence correspondence $\{(P, L): L\in (\bP^2)^{\ast}(\bF_q) \text{ and } P\in (L\cap C)(\bF_q)\}$. Since each $\bF_q$-point $P$ is contained in $q+1$ lines, 
$$
\#C(\bF_q) (q+1) = \sum_{P\in C(\bF_q)} (q+1) = \sum_{L} \# (L\cap C)(\bF_q)
$$
The sum on the right runs over all $q^2+q+1$ lines. Thus, an $\bF_q$-line on average contains
$$
\frac{\#C(\bF_q) (q+1)}{q^2+q+1} =  \frac{d(q-d+2)(q+1)}{q^2+q+1} > 
\frac{d(q-d+2)}{q+1} > d\left(1-\frac{d}{q+1}\right)
$$
$\bF_q$-points of $C$. As $q$ gets larger, this number approaches $d$. This heuristic suggests that Question \ref{Frob2} may have an affirmative answer.

\end{document}